\documentclass{amsart}

\usepackage{fancyhdr}
\usepackage{lastpage}
\usepackage{yhmath}

\newcommand{\bdis}{\begin{displaymath}}
\newcommand{\edis}{\end{displaymath}}
\newcommand{\be}{\begin{equation}}
\newcommand{\ee}{\end{equation}}

\newcommand{\mbb}{\mathbb}
\newcommand{\mcal}{\mathcal}

\newcommand{\vt}{\vartheta}

\newcommand{\btau}{\bar{\tau}}

\newcommand{\aZ}{\overset{*}{Z}}

\newcommand{\zf}{\zeta\left(\frac{1}{2}+it\right)}

\theoremstyle{definition}

\theoremstyle{remark}
\newtheorem{remark}[]{Remark}

\newtheorem*{mydef1}{{\bf Theorem}}

\newtheorem*{mydef5}{{\bf Lemma}}

\numberwithin{equation}{section}

\begin{document}

\title{Microscopic interpolation formula for the Riemann $Z(t)$-function and new algorithm for asymptotic solution of some Diophantine equation
with Leibnitz coefficients}

\author{Jan Moser}

\address{Department of Mathematical Analysis and Numerical Mathematics, Comenius University, Mlynska Dolina M105, 842 48 Bratislava, SLOVAKIA}

\email{jan.mozer@fmph.uniba.sk}


\begin{abstract}
In this paper we use the values of the Riemann $Z(t)$-function in order to construct certain quasi-orthonormal system of vectors. On this basis
we prove a formula for microscopic interpolation of the function $Z(t)$. Simultaneously we have obtained a new algorithm how to construct
asymptotic solutions of Diophantine equation with Leibnitz coefficients. This paper is English version of the paper \cite{3}
except the Appendix.
\end{abstract}

\maketitle

\section{Introduction}

Let (see \cite{4}, p. 79)
\be \label{1.1}
\begin{split}
 & Z(t)=e^{i\vt(t)}\zf, \\
 & \vt(t)=-\frac t2\ln\pi+\text{Im}\Gamma\left(\frac 14+i\frac t2\right),
\end{split}
\ee
and $\{ t_\nu\}$ be the sequence defined by the equation (see \cite{4}, p. 221)
\be \label{1.2}
\vt(t_\nu)=\pi\nu,\ \nu=1,2,\dots
\ee
In the paper \cite{2} we have obtained with the help of the Riemann-Siegel formula (see \cite{4}, p. 220)
\be \label{1.3}
Z(t)=2\sum_{n\leq\sqrt{\bar{t}}}\frac{1}{\sqrt{n}}\cos(\vt-t\ln n)+\mcal{O}(t^{-1/4}),\ \bar{t}=\frac{T}{2\pi},
\ee
the following formulae
\be \label{1.4}
\begin{split}
 & \sum_{T\leq t\nu\leq T+H}Z(t_\nu)Z(t_\nu+\bar{\tau}_k)\sim -\frac{1}{(4k+3)\pi^2}H\ln^2\frac{T}{2\pi}, \\
 & \sum_{T\leq t\nu\leq T+H}Z(t_\nu)Z(t_\nu+\bar{\bar{\tau}}_k)\sim \frac{1}{(4k+1)\pi^2}H\ln^2\frac{T}{2\pi},\ T\to\infty ,
\end{split}
\ee
where
\be \label{1.5}
\begin{split}
 & \bar{\tau}_k=\frac{(4k+3)\pi}{\ln\frac{T}{2\pi}}, \quad
 \bar{\bar{\tau}}_k=\frac{(4k+1)\pi}{\ln\frac{T}{2\pi}}, \\
 & k=0,1,\dots,K_0(T)=\mcal{O}(1),\ H=\sqrt{T}\ln T.
\end{split}
\ee

\begin{remark}
 The autocorrelative sum
 \be \label{1.6}
 \sum_{T\leq t\nu\leq T+H}Z(t_\nu)Z(t_\nu+\tau_k)
 \ee
 is oscillatory on the segment of the arithmetic sequence
 \bdis
 \bar{\tau}_0,\bar{\bar{\tau}}_0,\bar{\tau}_1,\bar{\bar{\tau}}_1,\dots,\bar{\tau}_{K_0},\bar{\bar{\tau}}_{K_0}.
 \edis
 as it follows by (\ref{1.4}).
\end{remark}

In this paper:
\begin{itemize}
 \item[(a)] we use the above mentioned oscillatory behavior of the sum (\ref{1.6}) in order to construc new kind
 (in the theory of the Riemann $Z(t)$-function) of quasi-orthonormal system of vectors;
 \item[(b)] we obtain, on basis of (a), a \emph{microscopic} interpolation formula for the function $Z(t)$ on a set
 of segments each of them has the length
 \be \label{1.7}
 < A\frac{\{\psi(T)\}^\epsilon}{\ln T}\to 0 \quad \text{as}\quad T\to\infty,
 \ee
 where $\psi(T)$ is an arbitrary fixed function of the type
 \bdis
 \ln_3T,\ln_4T,\dots;\quad \ln_3T=\ln\ln\ln T,\dots
 \edis
\end{itemize}

\section{Result}

\subsection{}

The following Theorem holds true.

\begin{mydef1}
 The formula
 \be \label{2.1}
 \begin{split}
  & \overset{*}{Z}(h_\nu)=\frac 4\pi\sum_{p=1}^L
  \frac{(-1)^{p+1}}{2p-1}\frac{\aZ(h_\nu+\bar{\tau}_p)+\aZ(h_\nu-\bar{\tau}_p)}{2}+ \\
  & + \mcal{O}\left\{ [\psi(T)]^{-\epsilon/4}\right\},\quad L=[(\psi(T))^\epsilon],
 \end{split}
\ee
where $[\dots]$ stands for the integer part, and
\be \label{2.2}
h_\nu=t_\nu+\frac{\pi}{\ln\frac{T}{2\pi}},\ \bar{\tau}_p=\frac{(2p-1)\pi}{\ln\frac{T}{2\pi}},\
\aZ(t)=\frac{Z(t)}{\sqrt{\ln\frac{T}{2\pi}}}
\ee
holds true for the values
\bdis
h_\nu\in [T,T+\bar{H}]
\edis
of the order
\be \label{2.3}
\sim \frac{1}{2\pi}\bar{H}\ln\frac{T}{2\pi},\quad T\to\infty;\ \bar{H}=\sqrt{T}\psi(T),
\ee
i. e. for \emph{the almost all} $h_\nu$.
\end{mydef1}

\begin{remark}
The result of our Theorem may be expressed, from the viewpoint of the theory of interpolation, as follows. If we assume the values
\bdis
\{\aZ(h_\nu\pm \bar{\tau}_p)\}_{p=1}^L
\edis
are given, then the formula (\ref{2.1}) expresses an approximation of the unknown value
\bdis
\aZ(h_\nu).
\edis
\end{remark}

\begin{remark}
For the discrete set of arguments in (\ref{2.1})
\bdis
\{ h_\nu-\bar{\tau}_L,h_\nu-\bar{\tau}_{L-1},\dots,h_\nu,h_\nu+\btau_1,\dots,h_\nu+\btau_L\}
\edis
we have (see (\ref{2.2}))
\be \label{2.4}
h_\nu+\btau_L-(h_\nu-\btau_L)=2\btau_L<A\frac{(\ln_3T)^\epsilon}{\ln T}\to 0,\ T\to\infty
\ee
(comp. (\ref{1.7})). Consequently, from (\ref{2.4}) the title \emph{microscopic} interpolation for (\ref{2.1}) follows.
\end{remark}

\subsection{}

The following Lemma is a basis for the proof of our Theorem.

\begin{mydef5}
If
\be \label{2.5}
\tau',\tau''=\mcal{O}\left(\frac{\psi^\epsilon}{\ln T}\right),
\ee
then we have the following formulae
\be \label{2.6}
\begin{split}
& \sum_{T\leq t_\nu\leq T+\bar{H}}Z(t_\nu+\tau')Z(t_\nu+\tau'')= \\
& = \frac{1}{2\pi}\frac{\sin\{ (\tau''-\tau')\ln P_0\}}{(\tau''-\tau')\ln P_0}\bar{H}\ln^2\frac{T}{2\pi}+\mcal{O}(\sqrt{T}\ln^2T),
\ \tau'\not=\tau'',
\end{split}
\ee
and
\be \label{2.7}
\begin{split}
& \sum_{T\leq t_\nu\leq T+\bar{H}}Z^2(t_\nu+\tau')=
\frac{1}{2\pi}\bar{H}\ln^2\frac{T}{2\pi}+\mcal{O}(\sqrt{T}\ln^2T),
\end{split}
\ee
where
\bdis
P_0=\sqrt{\frac{T}{2\pi}},
\edis
and these formulae are uniform in $\tau',\tau''$ (see the condition (\ref{2.5})).
\end{mydef5}

This Lemma follows from \cite{2}, (10), (11) in the case
\bdis
\begin{split}
& \tau\to\tau''-\tau',\ \tau=\mcal{O}\left(\frac{1}{\ln T}\right)\to\tau',\tau''=\mcal{O}\left(\frac{\psi^\epsilon}{\ln T}\right), \\
& H\to\bar{H}=\sqrt{T}\psi(T).
\end{split}
\edis

Next parts of this paper are ordered as follows. We define:
\begin{itemize}
\item[(a)] certain quasi-orthonormal system of vectors,
\item[(b)] an analogue of the Fourier coefficients for this case and, consequently, the asymptotic Fourier doefficients,
\item[(c)] an analogue of the trigonometric polynomial related with our vectors,
\item[(d)] corresponding mean square deviation.
\end{itemize}

After completion this program we prove the Theorem.

\section{Quasi-orthonormal system of vectors}

Since (see \cite{1}, (23))
\bdis
Q=Q(T,\bar{H})=\sum_{T\leq t_\nu\leq T+\bar{H}}1=\frac{1}{2\pi}\bar{H}\ln\frac{T}{2\pi}+\mcal{O}\left(\frac{\bar{H}^2}{T}\right),
\edis
then (see (\ref{2.3}))
\be \label{3.1}
Q\sim \frac{1}{2\pi}\bar{H}\ln\frac{T}{2\pi}; \ \bar{H}=\sqrt{T}\psi(T).
\ee
Next, we have in the case
\be \label{3.2}
\tau_p=\frac{2\pi}{\ln\frac{T}{2\pi}}p,\ p=-L+1,\dots,-1,0,1,\dots,L
\ee
that (see (\ref{2.5}), (\ref{2.6}), (\ref{3.2}))
\be \label{3.3}
\sum_{T\leq t_\nu\leq T+\bar{H}} Z^2(t_\nu+\tau_p)=\frac{1}{2\pi}\bar{H}\ln^2\frac{T}{2\pi}+\mcal{O}(\sqrt{T}\ln^2T),
\ee
and
\be \label{3.4}
\sum_{T\leq t_\nu\leq T+\bar{H}} Z(t_\nu+\tau_p)Z(t_\nu+\tau_p')=\mcal{O}(\sqrt{T}\ln^2T),\ p\not=p'.
\ee
Consequently, we obtain (see (\ref{2.2}), (\ref{3.1}) -- (\ref{3.4})) the following formula
\be \label{3.5}
\begin{split}
& \frac 1Q\sum_{T\leq t_\nu\leq T+\bar{H}} \aZ(t_\nu+\tau_p)\aZ(t_\nu+\tau_p')= \\
& = \left\{\begin{array}{rcl} \mcal{O}(\frac{1}{\psi}) & , & p\not=p', \\ 1+\mcal{O}(\frac{1}{\psi}) & , & p=p' . \end{array}\right.
\end{split}
\ee

\begin{remark}
We will call the property (\ref{3.5}) of the system of vectors
\be \label{3.6}
\{ \aZ(t_\nu+\tau_p)\},\ t_\nu\in [T,T+\bar{H}];\ p=-L+1,\dots,-1,0,1,\dots,L
\ee
as \emph{quasi-orthonormality}.
\end{remark}

\section{Analogue of Fourier coefficients and the classical Leibnitz series}

We define following numbers
\be \label{4.1}
\begin{split}
 & A_p=\frac 1Q\sum_{T\leq t_\nu\leq T+\bar{H}}\aZ(t_\nu+\tau^0)\aZ(t_\nu+\tau_p), \\
 & p=-L+1,\dots,-1,0,1,\dots,L
\end{split}
\ee
as an analogue of the Fourier coefficients of the vector
\be \label{4.2}
\aZ(t_\nu+\tau^0),\ T\leq t_\nu\leq T+\bar{H},\ \tau^0=\frac{\pi}{\ln\frac{T}{2\pi}}.
\ee
Since (see (\ref{3.2}), (\ref{4.2}))
\be \label{4.3}
\tau_p-\tau^0=\frac{(2p-1)\pi}{\ln\frac{T}{2\pi}}=\left(\pi p-\frac \pi2\right)\frac{1}{\ln P_0}; \
P_0=\sqrt{\frac{T}{2\pi}},
\ee
then we have (comp. (\ref{2.6}))
\be \label{4.4}
\frac{\sin\{ (\tau_p-\tau^0)\ln P_0\}}{(\tau_p-\tau^0)\ln P_0}=\frac \pi2\frac{(-1)^{p+1}}{2p-1}.
\ee
Consequently, we have from (\ref{2.6}) by (\ref{3.1}), (\ref{4.1}), (\ref{4.4}) that
\be \label{4.5}
A_p=\frac 2\pi\frac{(-1)^{p+1}}{2p-1}+\mcal{O}\left(\frac 1\psi\right).
\ee

\begin{remark}
 It is natural to call the numbers
 \be \label{4.6}
 \bar{A}_p=\frac 2\pi\frac{(-1)^{p+1}}{2p-1},\ p=-L+1,\dots,-1,0,1,\dots,L
 \ee
 as the asymptotic Fourier coefficients of the vector (\ref{2.4}).
\end{remark}

\begin{remark}
 Let us point out the presence of the members of the classical Leibnitz series
 \bdis
 \sum_{p=1}^\infty  \frac 2\pi\frac{(-1)^{p+1}}{2p-1} = \frac 12
 \edis
 in the notion of the asymptotic Fourier coefficients.
\end{remark}

\section{Analogue of the trigonometric polynomial}

Finally, we will define an analogue of the trigonometric polynomial $P_{2L}$, that corresponds to our quasi-orthonormal
system (\ref{3.6}) by the following way
\be \label{5.1}
P_{2L}[\aZ(t_\nu+\tau^0)]=\sum_{p=-L+1}^L \bar{A}_p\aZ(t_\nu+\tau_p).
\ee
Since from (\ref{4.6}) follows that
\be \label{5.2}
\bar{A}_p=\bar{A}_{-p+1},
\ee
then we have from (\ref{5.1}) by (\ref{4.6}), (\ref{5.2})
\be \label{5.3}
\begin{split}
 & P_{2L}=\sum_{p=1}^L \bar{A}_p\aZ(t_\nu+\tau_p)+\sum_{p=1}^L \bar{A}_{1-p}\aZ(t_\nu+\tau_{1-p})= \\
 & = \sum_{p=1}^L \bar{A}_p\{\aZ(t_\nu+\tau_p)+\aZ(t_\nu+\tau_{1-p})\} = \\
 & = \frac 2\pi\sum_{p=1}^L \frac{(-1)^{p+1}}{2p-1}\{\aZ(t_\nu+\tau_p)+\aZ(t_\nu+\tau_{1-p})\}.
\end{split}
\ee

\section{Mean square deviation and the classical Euler series}

\subsection{}

First of all, from the Euler series by (\ref{2.1}) we obtain
\be \label{6.1}
\begin{split}
 & \frac{\pi^2}{8}=\sum_{n=1}^\infty \frac{1}{(2p-1)^2}=\sum_{p=1}^L\frac{1}{(2p-1)^2}+\sum_{p=L+1}^\infty\frac{1}{(2p-1)^2}= \\
 & = \sum_{p=1}^L\frac{1}{(2p-1)^2}+\mcal{O}\left(\frac 1L\right)=\sum_{p=1}^L\frac{1}{(2p-1)^2}+\mcal{O}\left(\frac{1}{\psi^\epsilon}\right).
\end{split}
\ee
Since (see (\ref{4.6}))
\be \label{6.2}
\bar{A}_p^2=\frac{4}{\pi^2}\frac{1}{(2p-1)^2},
\ee
then we have (see (\ref{5.2}), (\ref{6.1}), (\ref{6.2}))
\be \label{6.3}
\begin{split}
 & \sum_{p=-L+1}^L\bar{A}_p^2=\frac{4}{\pi^2}\sum_{p=-L+1}^L\frac{1}{(2p-1)^2}=\frac{8}{\pi^2}\sum_{p=1}^L\frac{1}{(2p-1)^2}= \\
 & = \frac{8}{\pi^2}\left\{ \frac{\pi^2}{8}+\mcal{O}\left(\frac{1}{\psi^\epsilon}\right)\right\}=
 1+\mcal{O}\left(\frac{1}{\psi^\epsilon}\right),
\end{split}
\ee
and, of course,
\be \label{6.4}
\sum_{p,q=-L+1}^L |\bar{A}_p\bar{A}_q|=\mcal{O}(L^2)=\mcal{O}(\psi^{2\epsilon}).
\ee

\subsection{}

We define discrete mean square deviation $\Delta$ as follows
\be \label{6.5}
\Delta^2=\frac 1Q\sum_{T\leq t_\nu\leq T+\bar{H}}\{\aZ(t_nu+\tau^0)-P_{2L}\}^2.
\ee
Consequently, we have
\be \label{6.6}
\begin{split}
 & \Delta^2=\frac 1Q\sum_{(t_\nu)}\aZ^2(t_\nu+\tau^0)-\frac 2Q\sum_{(t_\nu)}\aZ(t_\nu+\tau^0)P_{2L}+\\
 & + \frac 1Q\sum_{(t_\nu)} P^2_{2L}=w_1+w_2+w_3.
\end{split}
\ee
From (\ref{2.7}), $\tau'=\tau''$, we obtain immediately (see (\ref{2.2}), (\ref{2.5}), (\ref{3.1}), (\ref{4.2})) that
\be \label{6.7}
w_1=1+\mcal{O}\left(\frac{1}{\psi}\right).
\ee

\subsection{}

Next, we obtain from (\ref{5.1}), (\ref{6.6}) by (\ref{4.1}), (\ref{4.5}), (\ref{4.6}), (\ref{6.3})
\be \label{6.8}
\begin{split}
 & w_2=-\frac 2Q\sum_{(t_\nu)}\aZ(t_\nu+\tau^0)\sum_{p=-L+1}^L \bar{A}_p\aZ(t_p+\tau_p)= \\
 & = -2\sum_{p=-L+1}^L \bar{A}_p\frac 1Q\sum_{(t_\nu)}\aZ(t_\nu+\tau^0)\aZ(t_\nu+\tau_p)= \\
 & = -2\sum_{p=-L+1}^L \bar{A}_pA_p=-2\sum_{p=-L+1}^L \bar{A}_p^2+\mcal{O}\left(\frac 1\psi\sum_{p=-L+1}^L 1\right)= \\
 & = -2+\mcal{O}\left(\frac{1}{\psi^\epsilon}\right)+\mcal{O}\left(\frac{1}{\psi^{1-\epsilon}}\right)= \\
 & = -2+\mcal{O}\left(\frac{1}{\psi^\epsilon}\right).
\end{split}
\ee

\subsection{}

Once again we have (see (\ref{3.5}), (\ref{5.1}), (\ref{6.3}), (\ref{6.4}))
\be \label{6.9}
\begin{split}
 & w_3=\sum_{p,q=-L+1}^L \bar{A}_p\bar{A}_q\frac 1Q\sum_{(t_\nu)}\aZ(t_\nu+\tau_p)\aZ(t_\nu+\tau_q)= \\
 & = \sum_{p=-L+1}^L \bar{A}_p^2\left\{ 1+\mcal{O}\left(\frac{1}{\psi}\right)\right\}+
 \mcal{O}\left\{\frac 1\psi\sum_{p,q=-L+1}^L |\bar{A}_p\bar{A}_q|\right\}= \\
 & = 1+\mcal{O}\left(\frac{1}{\psi^\epsilon}\right)+\mcal{O}\left(\frac{L^2}{\psi}\right)= \\
 & = 1+\mcal{O}\left(\frac{1}{\psi^\epsilon}\right).
\end{split}
\ee
Consequently, we obtain from (\ref{6.6}) by (\ref{6.7}) -- (\ref{6.9}) the following estimate
\be \label{6.10}
\Delta^2<\frac{A}{\{\psi(T)\}^\epsilon}.
\ee

\section{Finalisation of the proof of the Theorem}

Let $R(T)$ denote the number of points
\bdis
\bar{t}_\nu\in [T,T+\bar{H}]
\edis
that fulfill the inequality
\be \label{7.1}
|\aZ(\bar{t}_\nu+\tau^0)-P_{2L}[\aZ(\bar{t}_\nu+\tau^0)]|\geq \frac{1}{\psi^{\epsilon/4}}.
\ee
Our assertion is that
\be \label{7.2}
R(T)=o(\bar{H}\ln T),\quad T\to\infty.
\ee
Namely, if
\be \label{7.3}
R(T)>B\bar{H}\ln T
\ee
for every fixed $B>0$, then we have from (\ref{6.10}) by (\ref{6.5}), (\ref{7.1}), (\ref{7.3}) that
\bdis
\begin{split}
 & \frac{A}{\psi^\epsilon}\geq \frac 1Q\sum_{(\bar{t}_\nu)}\{\aZ(\bar{t}_\nu+\tau^0)-P_{2L}[\aZ(\bar{t}_\nu+\tau^0)]\}^2> \\
 & > \frac{B\bar{H}\ln T}{Q}\frac{1}{\psi^{\epsilon/2}}>\frac{C}{\psi^{\epsilon/2}},\ T\to\infty,
\end{split}
\edis
i. e. we have contradiction. Hence, our assertion (\ref{7.2}) is true. Consequently, if $Q_1(T)$ denotes the number of the points
\bdis
t_\nu\in [T,T+\bar{H}]
\edis
that fulfill the inequality
\be \label{7.4}
|\aZ(t_\nu+\tau^0)-P_{2L}[\aZ(t_\nu+\tau^0)]|<\frac{1}{\psi^{\epsilon/4}},
\ee
then (see (\ref{3.1}), (\ref{7.2}))
\be \label{7.5}
Q_1(T,\bar{H})\sim \frac{1}{2\pi}\bar{H}\ln\frac{T}{2\pi},\ T\to\infty,
\ee
i. e. the assertion of Theorem is true if we use the following notation
\bdis
t_\nu+\tau^0=h_\nu,\ \tau_p-\tau^0=\bar{\tau}_p,\ \tau_{1-p}-\tau^0=-\btau_p.
\edis

\appendix

\section{New algorithm for construction of an asymptotic solution to a Diophantine equation with Leibnitz coefficients}

First of all we express our formula (\ref{2.1}) in the following form
\be \label{A1}
\begin{split}
 & -\frac \pi4\aZ(h_\nu)+\sum_{p=1}^L \frac{(-1)^{p+1}}{2p-1}\frac{\aZ(h_\nu+\btau_p)+\aZ(h_\nu-\btau_p)}{2}= \\
 & = \mcal{O}(\psi^{-\epsilon/4})\sim 0,\ T\to\infty.
\end{split}
\ee
Next, we introduce (on this basis) the following Diophantine equation
\be \label{A2}
\sum_{p=0}^L a_p X_p=0,
\ee
where the numbers
\be \label{A3}
\begin{split}
 & a_0=-\frac \pi4,\ a_p=\frac{(-1)^{p+1}}{2p-1},\quad p=1,2,\dots,L; \\
 & a_0\in \mbb{R}\setminus\mbb{Q},\ a_p\in\mbb{Q}
\end{split}
\ee
are Leibnitz coefficients, i. e. the members of the classical Leibnitz series
\bdis
-\frac \pi4+\sum_{p=1}^\infty \frac{(-1)^{p+1}}{2p-1}=0.
\edis
Hence, from Theorem we have that for all sufficiently big
\bdis
T>0
\edis
there is unbounded set of asymptotic solutions
\be \label{A4}
\begin{split}
 & X_0=\aZ(h_\nu),\ X_p=\frac{\aZ(h_\nu+\btau_p)+\aZ(h_\nu-\btau_p)}{2}, \\
 & p=1,\dots,L,\ \left( L=[\{\psi\}^\epsilon]\to\infty \quad \text{as}\quad T\to\infty\right)
\end{split}
\ee
of Diophantine equation (\ref{A2}). Namely, for the number $Q_1(T)$ of
\bdis
h_\nu\in [T,T+\bar{H}]
\edis
that generate asymptotic solutions, we have (see (\ref{2.3}), (\ref{7.5}))
\bdis
Q_1(T)\sim\frac{1}{2\pi}\bar{H}\ln\frac{T}{2\pi}=\frac{1}{2\pi}\sqrt{T}\psi(T)\ln\frac{T}{2\pi}\to\infty \quad \text{as} \quad
T\to\infty.
\edis

\begin{remark}
 It is clear, that from the point of view of this Appendix, the content of this paper may be interpreted as new algorithm
 how to construct certain set of asymptotic solution (\ref{A4}) of the Diophantine equation (\ref{A2}) with the
 Leibnitz coefficients (\ref{A3}). In order to attain this result we have used the values of the Riemann
 $Z(t)$-function considering its argument $t$ belongs to the \emph{microscopic} segment (see (\ref{2.4})).
\end{remark}

\thanks{I would like to thank Michal Demetrian for helping me with the electronic version of this work.}

\end{document}